\documentclass[11pt]{amsart}

\usepackage[T1]{fontenc}
\usepackage[utf8]{inputenc}
\usepackage{lmodern}
\usepackage{microtype}

\usepackage{CJKutf8}
\usepackage{mathtools,amssymb,amsfonts,amsthm}
\usepackage{euscript}

\usepackage{graphicx}
\usepackage{tikz}
\usepackage{tikz-cd}
\usepackage[arrow,matrix,curve]{xy}

\usepackage[
backend=biber,
style=alphabetic,
maxalphanames=10,
maxbibnames=10,
sorting=anyt,
]{biblatex}

\usepackage[normalem]{ulem} 
\usepackage{xcolor}

\usepackage[margin=1in]{geometry}
\usepackage{setspace}
\makeatletter
\g@addto@macro\normalsize{%
  \setlength\abovedisplayskip{12pt}%
  \setlength\belowdisplayskip{12pt}%
  \setlength\abovedisplayshortskip{10pt}%
  \setlength\belowdisplayshortskip{10pt}%
}
\makeatother

\allowdisplaybreaks[2]

\DeclareMathOperator{\T}{T}

\newcommand{\MM}{\mathcal{M}} 

\newcommand{\bH}{{\mathbb{H}}}
\newcommand{\bR}{{\mathbb{R}}}

\newcommand{\bC}{{\mathbb{C}}}

\newcommand\al\alpha
\newcommand\gam\gamma
\newcommand\del\delta
\newcommand\kap\kappa
\newcommand\lam\lambda
\newcommand\sig\sigma
\newcommand\veps\varepsilon
\newcommand\Gam\varGamma
\newcommand\bra\langle
\newcommand\ket\rangle

\numberwithin{equation}{section}
\newtheorem{thm}{Theorem}[section]

\newtheorem{cor}[thm]{Corollary}
\newtheorem{pro}[thm]{Proposition}
\theoremstyle{definition}
\newtheorem{dfn}[thm]{Definition}

\newcommand\Proof{\vspace{-10pt}\noindent {\bf Proof.} }

\tikzcdset{column sep=large, row sep=large}

\begin{document}
\title{Anomaly Reparametrization of the Ligon--Schaaf Regularization in the Kepler problem}
\author{Li-Chun Hsu}
\maketitle

\begin{abstract}
We revisit the Ligon--Schaaf regularization of the Kepler problem and identify the geometric origin of the rotation appearing in their transformation.
We show that this rotation is determined by the eccentric anomaly of the Kepler motion,
providing a transparent dynamical interpretation of the angle that renders the Kepler flow
uniform on $T^{*}S^{3}$. 
Building on this insight, we extend the construction to positive and zero energies
via the corresponding hyperbolic and parabolic anomalies, obtaining a unified geometric
description of the Kepler flow across all energy levels.
\end{abstract}
\tableofcontents

\clearpage

\section{Introduction}
The classical Kepler problem
\[
H(q,p)=\tfrac{1}{2}|p|^2-\mu |q|^{-1},
\qquad (q,p)\in T\bR^3\cong \bR^3\times \bR^3,
\]
admits several conserved quantities.
The rotational invariance of the Hamiltonian implies conservation of the angular momentum
$$L=q\times p.$$
Beyond this manifest symmetry, the system possesses an additional hidden symmetry encoded by the Laplace–Runge–Lenz vector
$$A=p\times L-\frac{q}{|q|}.$$
This conserved quantity can be traced back to the work of Jacob Hermann and his teacher Johann Bernoulli in 1710.
On the negative-energy region, this enhanced symmetry allows the dynamics to be regularized
and identified with geodesic motion on a compact manifold.

Moser's celebrated construction \cite{Moser1970Kepler} realizes each fixed negative energy level
$H^{-1}(c)$, $c<0$, as the unit cotangent bundle $T^*_1S^3$ of the $3$--sphere,
transforming Keplerian motion into the geodesic flow of $S^{3}$.
However, the Moser regularization is inherently \emph{level-wise}:
the symplectomorphism is defined separately on each fixed negative energy surface.

To obtain a global correspondence valid for all negative energies at once,
Ligon and Schaaf \cite{LigonSchaaf1976Kepler} introduced an additional rotation into Moser construction.
This yields a global symplectomorphism from the whole negative-energy region to
an open subset of $T^*S^3\times\mathbb{R}$ such that the Kepler flow becomes a
reparametrized geodesic flow. While their regularization map has been analyzed
from algebraic and symplectic viewpoints, the geometric origin and physical
meaning of this extra ``rotation'' have remained somewhat unclear in the literature.
In particular, Ligon and Schaaf briefly remark that their angle is related to the
Kepler equation, but the precise interpretation
has not been made fully explicit.

The first goal of this paper is to clarify this point.
We show that the rotation introduced by Ligon and Schaaf 
is precisely a reparametrization of the Kepler flow by the \emph{eccentric anomaly}.
More concretely, starting from Moser map and following the evolution of a Kepler orbit,
we prove that the additional angle appearing in the Ligon--Schaaf transformation is exactly
\[
    \varepsilon - M,
\]
where $\varepsilon$ is the eccentric anomaly and $M$ is the mean anomaly.
Consequently, the Ligon--Schaaf map is not merely a technical modification that enforces
energy--uniformity: it is a geometrically natural \emph{anomaly--based reparametrization}
that converts the non--uniform Kepler motion into a uniform geodesic rotation on~$S^{3}$.
This makes precise, in a fully symplectic setting, the observation (implicit in the
original Ligon--Schaaf article \cite{LigonSchaaf1976Kepler}) that their auxiliary angle measures the deviation
of the eccentric from the mean anomaly, or equivalently, that their equation is
a generalized Kepler equation.

A second contribution of this work is to extend this interpretation to the
positive and zero energy cases.
For positive energy, we revisit Belbruno's hyperbolic regularization and incorporate
a Ligon--Schaaf--type rotation, showing that the resulting map takes the Kepler flow to
the geodesic flow on hyperbolic space $H^{3}$, with the rotation angle given by the
\emph{hyperbolic anomaly}.
For zero energy, we treat the Euclidean limit and obtain an analogous formulation
based on the parabolic anomaly.
Thus, all energy levels of the Kepler problem---elliptic, hyperbolic, and parabolic---fit
into a unified geometric framework governed by a single principle:
\[
    \text{Moser regularization} 
    \quad+\quad
    \text{anomaly reparametrization}
   =\text{geometric uniformization of Kepler flow}.
\]

The constructions of Moser, Ligon--Schaaf, and Belbruno are individually known.
Our contribution is to extract and isolate the common geometric mechanism and show 
that it is naturally expressed in terms of the anomaly parameters of celestial mechanics.

\medskip
The rest of the paper is organized as follows.
In Section~\ref{sec:scaling}, we recall the scaling symmetry of the Kepler Hamiltonian
and the auxiliary Hamiltonians $K_\pm$ and $K_0$ that normalize the different energy levels.
Section~\ref{sec:negative} reviews the Moser regularization for negative energies, expressed
in quaternionic language \cite{Knorrer1997Kepler}, and the induced symmetry and moment map.
In Section~\ref{sec:LS-negative} we introduce the Ligon--Schaaf rotation, identify its
angle with the eccentric anomaly and show how the Kepler flow becomes uniform
on $T_1^*S^3$. Section~\ref{sec:positive} is devoted to positive energy:
we recall Belbruno's hyperbolic regularization, extend it to a global Ligon--Schaaf--type map,
and identify the role of the hyperbolic anomaly.
Finally, in Section~\ref{sec:zero} we treat the parabolic case using a
Euclidean regularization and discuss the corresponding anomaly interpretation.

\clearpage

\section{Scaling and normalized Hamiltonians}
\label{sec:scaling}

Let $T^*\bR^3$ be the cotangent bundle of $\bR^3$. We work with the Kepler Hamiltonian
\[
H(q,p)=\tfrac{1}{2}|p|^2-\mu |q|^{-1},\qquad (q,p)\in T\bR^3\cong \bR^3\times \bR^3.
\]

The Hamiltonian equations for $H$ are invariant under the scaling
\[
\rho\cdot(t,q,p):=(\rho^3 t,\rho^2 q,\rho^{-1}p),\qquad \rho>0,
\]
under which $H$ transforms as $H\mapsto \rho^{-2}H$.
Thus, to understand all orbits up to scaling, it is enough to study the level sets
\[
H^{-1}\!\left(-\tfrac12\right),\qquad
H^{-1}(0),\qquad
H^{-1}\!\left(\tfrac12\right).
\]

It is convenient to introduce the following auxiliary Hamiltonians,
obtained by multiplying $H-c$ by $2|q|$ for $c\in\{0,\pm\tfrac12\}$:
\begin{align*}
K_\pm(q,p)&=2|q|\big(H(q,p)\mp\tfrac12\big)
=|q|(|p|^2\pm 1)-2\mu,\\
K_0(q,p)&=2|q|\big(H(q,p)-0\big)
=|q|\,|p|^2-2\mu.
\end{align*}
Without loss of generality, we set $\mu=1$ by choosing suitable units.
Then
\[
K_+^{-1}(0)=H^{-1}\!\left(-\tfrac{1}{2}\right),\qquad
K_-^{-1}(0)=H^{-1}\!\left(\tfrac{1}{2}\right),\qquad
K_0^{-1}(0)=H^{-1}(0).
\]

The normalized Hamiltonians $K_\pm$ and $K_0$ play a central role:
on the corresponding energy surfaces $K_\pm^{-1}(0)$ and $K_0^{-1}(0)$, their Hamiltonian vector fields are proportional to $X_H$, where $X_H$ is the Hamiltonian vector field of $H$, by a simple factor
depending only on $|q|$; see Propositions~\ref{pro:Kplus-XH} and~\ref{pro:Kminus-XH}
below.

\clearpage

\section{Anomalies}

\noindent\textbf{Notation.}
Throughout this section, we consider planar Keplerian motion under a central inverse--square force.
The position of the orbiting body is described in polar coordinates $(r,\theta)$ with respect to the focus of the conic, which is placed at the origin.
The eccentricity is denoted by $e\ge 0$, and $\mathfrak{p}>0$ denotes the semi-latus rectum.
For elliptic orbits ($0<e<1$), $a>0$ denotes the semi-major axis and $\mathfrak{p}=a(1-e^2)$.
The constants $\alpha$ and $\mu$ denote the strength of the central potential and the reduced mass, respectively.
Time is measured from the periapsis passage and denoted by $t-t_p$.
Angular motion is assumed to be monotone, i.e.,$\dot{\theta}>0$. For completeness, we briefly review the notion of anomaly, mainly following
Arnold, Kozlov, and Neishtadt \cite{Arnold1997MathematicalAO}.

\medskip
Once the geometric shape of a Keplerian orbit is fixed by its orbital elements, such as the semi-major axis and the eccentricity, it remains to parametrize the position of the body along the orbit as a function of time.
For this purpose, it is convenient to introduce several auxiliary angular parameters, collectively referred to as \emph{anomalies}.
These parameters lead to simple expressions for the radial distance and for the swept area, and hence provide natural coordinates for describing orbital motion.

Throughout this section, we fix the periapsis direction by setting the phase angle $\theta_0=0$.
The Keplerian orbit is then described by
\[
r(\theta)=\frac{\mathfrak{p}}{1+e\cos\theta}.
\]
This $\theta$ is called the true anomaly.

We choose the coordinate axes x and y along the major axes of the conic section representing the orbit. 
One convenient way of writing the conic equation is 
\begin{align*}
    e=1&\Rightarrow &\!\!\!\!\!\!\!\!\!\!\!\!\!\!\!\!\!\!\!\!\!\!\!\!\!\!\!\!\!\!\!\!\!\!\!\!\!\!\!\!\!\!\!\!\!\!\!\!\!\!\!\!\!\!\!\!&y^2+2\mathfrak{p}x=\mathfrak{p}^2,\\
    e\neq1&\Rightarrow &\!\!\!\!\!\!\!\!\!\!\!\!\!\!\!\!\!\!\!\!\!\!\!\!\!\!\!\!\!\!\!\!\!\!\!\!\!\!\!\!\!\!\!\!\!\!\!\!\!\!\!\!\!\!\!\!&\frac{(x+ae)^2}{a^2}\pm\frac{y^2}{b^2}=1
\end{align*}
The equation of the orbit can be represented in the following parametric form:
\begin{align*}
&x=a(\cos\psi-e),& &y=a\sqrt{1-e^2}\sin\psi& &(0\leq e <1) &if\ h<0;\\
&x=a(\cosh\psi_h-e),& &y=a\sqrt{e^2-1}\sinh\psi_h&  &(e>1) &if\ h>0;\\
&x=\frac{p}{2}(1-\psi_p^2),& &y=p\psi_p&  &(e=1) &if\ h=0.
\end{align*}
These angles are called eccentric anomaly, hyperbolic anomaly, and parabolic anomaly.

And by simple computation, we can get the following equations:
\begin{align*}
    &(0\leq e<1) &\!\!\!\!\!\!\!\!\!\!\!\!\!\!\!\!\!\!\!\!\!\!\!\!\!\!\!\!\!\!\!\!\!\!&\tan\frac{\theta}{2}=\sqrt{\frac{1+e}{1-e}}\tan\frac{\psi}{2}\\
    &(e>1) &\!\!\!\!\!\!\!\!\!\!\!\!\!\!\!\!\!\!\!\!\!\!\!\!\!\!\!\!\!\!\!\!\!\!&\tan\frac{\theta}{2}=\sqrt{\frac{1+e}{1-e}}\tanh\frac{\psi_h}{2}\\
    &(e=1) &\!\!\!\!\!\!\!\!\!\!\!\!\!\!\!\!\!\!\!\!\!\!\!\!\!\!\!\!\!\!\!\!\!\!&\tan\frac{\theta}{2}=\psi_p
\end{align*}

In the ellipse case, let $\omega$ be the motion, i.e., the average angular velocity $\omega:=\frac{2\pi}{T}$, where $T$ is the period. The mean anomaly is defined by $$M=\omega(t-t_p)$$
The swept-area computation shows that $M$ is naturally paired with $\psi$ through the normal form
\[
M=\psi-e\sin\psi,
\]
which is known as the Kepler Equation.

However, in the hyperbolic and parabolic cases the motion is not periodic, so the mean anomaly should no longer be interpreted as an angular variable, but rather as a time–area parameter. According to Kepler's third law, we may rewrite the mean motion as $$\omega=\sqrt{\frac{\alpha}{\mu a^3}}$$
which doesn't include the period, and inspires us to define mean anomaly in the other two cases.

\medskip
\noindent\textbf{Hyperbolic anomaly.}
We define the \emph{hyperbolic mean motion} by
\[
\omega_h:=\sqrt{\frac{\alpha}{\mu(-a)^3}},
\]
and we define the \emph{hyperbolic mean anomaly} by
\[
M_h:=\omega_h\,(t-t_p)
\]
Also, the swept-area computation shows that $M_h$ is naturally paired with $\psi_h$ through the normal form $$M_h=e\sinh\psi_h-\psi_h$$
This expression is the hyperbolic Kepler equation.

\medskip
\noindent\textbf{Parabolic anomaly.}
Similarily, we define the parabolic mean motion and parabolic mean anomaly
\[
\omega_p:=2\sqrt{\frac{\alpha}{\mu\,\mathfrak{p}^3}},
\qquad
M_p:=\omega_p\,(t-t_p),
\]

Kepler's second law again implies that the elapsed time from periapsis passage is proportional to the swept area.
This leads naturally to the polynomial normal form
\[
M_p =\psi_p+\frac{1}{3}\psi_p^3.
\]

\medskip
In summary, starting from the elliptic parametrization and Kepler's second law, one is naturally led to the canonical normal forms
\[
M=\psi-e\sin\psi,\qquad
M_h=e\sinh\psi_h-\psi_h,\qquad
M_p=\psi_p+\frac{1}{3}\psi_p^3,
\]
which we adopt as the defining relations for the mean anomalies in the elliptic, hyperbolic, and parabolic cases, respectively.
These anomalies will be used in the sequel as convenient coordinates for describing Keplerian motion.

\clearpage

\section{Negative energy: Moser and Ligon--Schaaf}
\label{sec:negative}

\subsection{Moser regularization}
\label{subsec:Moser}
In this subsection, we review Moser regularization introduced in \cite{Moser1970Kepler}.
We begin with the stereographic projection. Let
\[
S^3_{np}:=S^3\setminus\{(1,0,0,0)\}
\]
be the $3$--sphere minus the north pole. Define
\[
s\colon S^3_{np}\to\bR^3,\qquad s(x)=\frac{1}{1-x_0}(x_1,x_2,x_3).
\]
This is a diffeomorphism. Any diffeomorphism $f\colon M_1\to M_2$ induces a natural
symplectomorphism on cotangent bundles,
\[
f_\#=(f,((df)^*)^{-1})\colon T^*M_1\stackrel{\rm sympl.}{\longrightarrow}T^*M_2,
\]
hence $s_\#$ is symplectic.

Let
\[
g\colon T^*\bR^3\to T^*\bR^3,\qquad g(q,p)=(-p,q).
\]
Then $g$ is symplectic, and so is the composition
\[
g\circ s_\#\colon T^*S^3_{np}\to T^*\bR^3.
\]
In coordinates, one checks that
\[
g\circ s_\#(x,y)
=\Big(-(1-x_0)(y_1,y_2,y_3)-y_0(x_1,x_2,x_3),\ (1-x_0)^{-1}(x_1,x_2,x_3)\Big).
\]
Writing $g\circ s_\#(x,y)=(q,p)$, the \emph{Moser regularization} is the restriction
\[
\MM:=g\circ s_\#\big|_{T_1^*S^3_{np}}\colon T_1^*S^3_{np}\to T^*\bR^3.
\]
where $T_1^*S^3_{np}:=\{(x,y)\in T^*S^3_{np}\mid x\in S_{np}^3, y\in T^*_xS_{np}^3, |y|=1\}$
\begin{pro}[\cite{Moser1970Kepler}]
\label{pro:Moser-inverse}
The map $\MM$ is a diffeomorphism onto $K_+^{-1}(0)=H^{-1}(-\tfrac12)$,
and its inverse is given by
\[
\MM^{-1}(q,p)=\Big[\Big(\frac{|p|^2-1}{|p|^2+1},\frac{2p}{|p|^2+1}\Big),
\Big(-\bra q,p\ket,-\tfrac{|p|^2+1}{2}\,q+\bra q,p\ket\, p\Big)\Big].
\]
\end{pro}

\Proof
A direct computation shows that for $(x,y)\in T_1^*S^3_{np}$ one has
\[
|q|=|y|\,(1-x_0),\qquad |p|^2+1=2(1-x_0)^{-1}.
\]
Hence
\[
K_+(q,p)=|q|(|p|^2+1)-2=2|y|-2,
\]
so that $K_+(q,p)=0$ is equivalent to $|y|=1$, i.e.\ $(x,y)\in T_1^*S^3_{np}$.
The stated formula for $\MM^{-1}$ is obtained by solving the above relations
for $(x,y)$ in terms of $(q,p)$ and verifying that the result lies in $T_1^*S^3$.
\qed

Thus, Moser regularization identifies the fixed negative-energy surface
$H^{-1}(-\tfrac12)$ with the unit cotangent bundle $T_1^*S^3$,
and the Kepler flow with the geodesic flow on $S^3$ for a suitable time change.
We now rewrite this construction in quaternionic language, which is particularly
convenient for describing symmetry.

\subsection{Quaternionic description and symmetry}

We identify $\bR^4$ with the quaternions $\bH$ via
\[
(x_0,x_1,x_2,x_3)\longleftrightarrow x_0+ix_1+jx_2+kx_3,
\]
and let $\bH_{\mathrm{pure}}:=\{x\in\bH\mid x=-\bar x\}$, 
$\bH_{\mathrm{unit}}:=\{x\in\bH\mid x\bar x=1\}$
Under this identification, we can regard
\[
T_1^*S^3
=\{(x,y)\in\bH^2: |x|=1,\ |y|=1,\ xy^*+yx^*=0\},
\]
and the Moser map becomes
\[
\MM\colon T_1^*S^3_{np}\to\bH_{\mathrm{pure}}^2,\qquad
(x,y)\mapsto\Big(\tfrac{1}{2}(y-y^*)+\tfrac{1}{2}(x^*y-xy^*),\ \frac{1+x}{1-x}\Big).
\]

Let $G_+:=\bH_{\mathrm{unit}}^2$ act on $T^*S^3$ by
\[
(\alpha_1,\alpha_2)\cdot(x,y)=(\alpha_1x\alpha_2^*,\alpha_1y\alpha_2^*).
\]
By  the explicit formula above, this induces an
action $*$ of $G_+$ on $\bH_{\mathrm{pure}}^2\cong T^*\bR^3$, characterized by the property
\[
\MM\big((\alpha_1,\alpha_2)\cdot(x,y)\big)
=(\alpha_1,\alpha_2)*\MM(x,y).
\]

\begin{thm}[\cite{Knorrer1997Kepler}, Theorem 3.1]
\label{thm:Gplus-action}
The above $G_+$ action is explicitly given by
\[
(\alpha_1,\alpha_2)*(q,p)
=\big((bp+a)q(bp+a)^*,\ (ap+b)(bp+a)^{-1}\big),
\]
where $a=\tfrac{\alpha_1+\alpha_2}{2}$ and $b=\tfrac{\alpha_2-\alpha_1}{2}$.
The action is symplectic $G_+$–action on $T^*\bR^3\cong\bH_{\mathrm{pure}}^2$ and preserves the Hamiltonian $K_+$.
\end{thm}

Note that $(2i,0)\in K_+^{-1}(0)$. If $(\alpha_1,\alpha_2)*(2i,0)=(2i,0)$, then
$(2aia^*,ba^{-1})=(2i,0)$, hence $b=0$ and $aia^*=i$, so $|a|=1$ with $a\in\bC$.
Thus $\dim(\mathrm{stab}_{G_+}(2i,0))=\{(\alpha,\alpha)\mid\alpha\in\bC,|\alpha|=1\}=1$.

Define $\phi_+\colon G_+\to K_+^{-1}(0)$ by $\phi_+(\alpha_1,\alpha_2)=(\alpha_1,\alpha_2)*(2i,0)$.
It is surjective, since $\dim K_+^{-1}(0)=5$. Because $\MM(-1,i)=(2i,0)$ and $\MM$ is a
diffeomorphism $T_1^*S^3_{np}\to K_+^{-1}(0)$, we have
\[
\MM(-\alpha_1\alpha_2^*,\alpha_1 i\alpha_2^*)=\phi_+(\alpha_1,\alpha_2),
\qquad
\phi_+(G_+)=K_+^{-1}(0)\ \Longrightarrow\ 
\{(\alpha_1,\alpha_2)\cdot(-1,i)\}=T_1^*S^3_{np}.
\]
Let $I_+:=\{(u,v)\in\bH\setminus\{0\}\times\bH\mid|u|^2+|v|^2=1,u^*v+v^*u=0\}$. \\
Composing $g_+\!:I_+\rightarrow G_+,\ (u,v)\mapsto(u-v,u+v)$ with $\phi_+$ recovers the
Kustaanheimo--Stiefel map \cite{KustaanheimoStiefel1965KS}. We obtain the diagram
\[
\begin{tikzcd}
I_+ \arrow[rd, "KS"]\arrow[d, "g_+"]&        & \mathbb{R} &        \\
G_+ \arrow[r, "\phi_+"]  & K_+^{-1}(0) \arrow[ru, "K_+"] \arrow[rr, "\MM"] 
\arrow["G_+"', loop, distance=2em, in=305, out=235] &            & 
T_1^*S_{np}^3 \arrow[lu, "2|y|-2"'] \arrow["G_+"', loop, distance=2em, in=305, out=235]
\end{tikzcd}
\]

\subsection{Momentum map and integrals of motion}

By Noether's theorem, the $G_+$–action yields a conserved momentum map.
We compute it explicitly and identify it with the Laplace--Runge--Lenz vector $A=(A_1, A_2, A_3)=p\times L-\frac{q}{|q|}$ and
the angular momentum $L=(L_1, L_2 ,L_3)=q\times p$.

\begin{pro}
\label{pro:moment-map-negative}
The Lie algebra $\mathfrak{g}_+$ of the Lie group $G_+$ is
\[
\mathfrak{g}_+=\bH_{\mathrm{pure}}\oplus\bH_{\mathrm{pure}},
\]
and the vector field at $(x,y)\in T^*S^3$ associated to $(a,b)\in\mathfrak{g}_+$
is
\[
V^{(a,b)}(x,y)=(ax-xb,\ ay-yb).
\]
The corresponding momentum map $\mu\colon T^*S^3\to\mathfrak{g}_+^*$ satisfies
\[
\mu_{(a,b)}(x,y)
=\bra y,ax-xb\ket=(a+b)\cdot (\varepsilon_{ijk}x_jy_k)+(b-a)\cdot (x_iy_0-x_0y_i),
\]
\end{pro}

\Proof
The description of $\mathfrak{g}_+$ follows from the condition
\[
(1+\varepsilon a)(1+\varepsilon a)^*=1
\]
to first order in $\varepsilon$, which implies $a+\bar a=0$, and similarly for $b$.
The infinitesimal generator is obtained by differentiating the group action
at $\varepsilon=0$.
The expression for $\mu$ is the contraction of $V^{(a,b)}$ with the canonical $1$--form.
The identification with $A$ and $L$ on the Kepler energy surface
$H=-\tfrac12$ is obtained by composing with $\MM$ and using the explicit
formulas in Theorem~\ref{thm:AL-Moser} below.
\qed

\begin{thm}\cite{Moser1970Kepler}
\label{thm:AL-Moser}
\[
L_i\circ \MM=\varepsilon_{ijk}x_jy_k:= L'_i,\qquad
[A_i-(H+\tfrac{1}{2})q_i]\circ \MM=x_iy_0-x_0y_i:=A'_i,\qquad
K_+\circ \MM=2|y|-2.
\]
\end{thm}

\Proof
We have already computed $K_+\circ\MM$ above.
For the angular momentum,
\begin{align*}
L_i\circ \MM&=(\varepsilon_{ijk}q_jp_k)\circ \MM\\
&=-\varepsilon_{ijk}((1-x_0)y_i+y_0x_j)(1-x_0)^{-1}x_k\\
&=-\varepsilon_{ijk}y_jx_k=\varepsilon_{ijk}x_jy_k\\[6pt]
\end{align*}

For the Laplace--Runge--Lenz vector $A$,

\begin{align*}
A_i\circ \MM&=(q_i|p|^2-p_i\bra q,p\ket-\tfrac{q_i}{|q|})\circ \MM\\
&=-[(1-x_0)y_i+y_0x_i]\cdot\frac{1+x_0}{1-x_0}+\frac{y_0}{1-x_0}x_i-\frac{-(1-x_0)y_i-y_0x_i}{|y|(1-x_0)}\\
&=(\tfrac{1}{|y|}-x_0-1)y_i+\frac{y_0x_i}{1-x_0}(\tfrac{1}{|y|}-x_0).\\[6pt]
(H+\tfrac{1}{2})q_i\circ \MM&=(\tfrac{-1}{1-x_0}(1-\tfrac{1}{|y|})[(1-x_0)y_i+y_0x_i])\circ \MM\\
&=-(1-\tfrac{1}{|y|})y_i-\tfrac{1}{1-x_0}(1-\tfrac{1}{|y|})y_0x_i
\end{align*}
\[\Rightarrow(A_i-(H+\tfrac{1}{2})q_i)\circ \MM=x_iy_0-x_0y_i.\]\qed

\subsection{Poisson bracket and time reparametrization}

In this subsection, we investigate the relation between the Hamiltonian flows of $H$ and $K_+$.

\begin{pro}
\label{pro:Kplus-XH}
For the Hamiltonians
\[
H(q,p) = \tfrac{1}{2}|p|^2 - \tfrac{1}{|q|}, 
\qquad 
K_+(q,p) = |q|(|p|^2 + 1) - 2,
\]
one has
\[
\{K_+, H\} = \frac{\bra q,p\ket}{|q|^2}\, K_+(q,p).
\]
In particular, on $K_+^{-1}(0)=H^{-1}(-\tfrac12)$ one has $\{K_+,H\}=0$, and
\[
X_{K_+}=2|q|\,X_H\quad\text{on }K_+^{-1}(0).
\]
\end{pro}

\Proof
We compute
\[
\frac{\partial K_+}{\partial q_i} = (|p|^2 + 1)\frac{q_i}{|q|},
\qquad 
\frac{\partial K_+}{\partial p_i} = 2|q|\,p_i,
\]
\[
\frac{\partial H}{\partial p_i} = p_i,
\qquad
\frac{\partial H}{\partial q_i} = \frac{q_i}{|q|^3}.
\]
Thus
\begin{align*}
\{K_+, H\} 
&= \sum_i 
\left(
  \frac{\partial K_+}{\partial q_i}\frac{\partial H}{\partial p_i}
  - 
  \frac{\partial K_+}{\partial p_i}\frac{\partial H}{\partial q_i}
\right)\\
&=(|p|^2 + 1)\frac{\bra q,p\ket}{|q|}
- 2|q|\frac{\bra q,p\ket}{|q|^3} \\
&=
\frac{\bra q,p\ket}{|q|^2}\big[(|p|^2 + 1)|q| - 2\big]
=\frac{\bra q,p\ket}{|q|^2}K_+(q,p).
\end{align*}
This vanishes on $K_+^{-1}(0)$.
A direct computation gives
\[
X_H= 
\left(p, -\frac{q}{|q|^3}\right),\qquad
X_{K_+} = 
\left(2|q|p, -\frac{|p|^2 + 1}{|q|}q\right).
\]
On $K_+^{-1}(0)$ we have
\[
\frac{|p|^2+1}{|q|}=\frac{2}{|q|^2},
\]
so
\[
X_{K_+}=2|q|\Big(p,-\frac{q}{|q|^3}\Big)=2|q|\,X_H.
\]
\qed

Thus $K_+$ generates the same trajectories as $H$ on $H=-\tfrac12$,
but parametrized by a rescaled time variable.

\subsection{Energy-uniform Moser map}

Moser original construction \cite{Moser1970Kepler} is defined only on the fixed energy surface
$H^{-1}(-\tfrac12)$. To pass to an energy-uniform version, we adapt the Moser map
using the scaling symmetry.

Let
\[
\Sigma_-:=\{(q,p)\in T^*(\bR^3\setminus\{0\})\mid H(q,p)<0\}
\]
be the negative-energy region, and
\[
N_-:=\{(x,y)\in (T^*S^3)^\times\mid x_0\neq 1\}
\]
be the corresponding subset in $T^*S^3$.
Define
\[S\colon\Sigma_-\to \T^*_1S^3\times\bR_{>0}\]
\[
S(q,p)=\Big[(r_0,\vec r),\ (s_0,\vec s),\ \sqrt{-2H}\Big]
\]
with the definitions
\[
r_0:=|p|^2|q|-1,\qquad \vec r:=\sqrt{-2H}\,|q|p,
\]
\[
s_0:=-\sqrt{-2H}\,\bra q,p\ket,\qquad
\vec s:= -\Big(\frac{q}{|q|}-\bra q,p\ket\,p\Big).
\]

\begin{pro}
When $H=-\tfrac{1}{2}$, the map $S$ reduces to
\[
S|_{H=-1/2}=(\MM^{-1},\;1),
\]
where we regard $S$ as taking values in $T_1^*S^3\times\bR$. Moreover, $S$ is the unique map with this property and is invariant under the scaling in Chapter 2.
\end{pro}

\Proof
If $H=-\frac{1}{2}$ then
\[
\frac{1}{2}|p|^2-|q|^{-1}=-\frac{1}{2}
\quad\Longrightarrow\quad
|p|^2-\frac{2}{|q|}=-1 
\quad\Longrightarrow\quad
|q|=\frac{2}{|p|^2+1}.
\]
Substituting into the definitions of $r_0,\bar r,s_0,\bar s$ gives
\[
S(q,p)
=\Big[(\tfrac{|p|^2-1}{|p|^2+1},\tfrac{2p}{|p|^2+1}),
(-\bra q,p\ket,-\tfrac{|p|^2+1}{2}q+\bra q,p\ket \, p),1\Big],
\]
which is precisely $\MM(q,p)$ together with the value $1$ of $\sqrt{-2H}$.
\qed

One checks that $S$ defines a symplectomorphism from $N_-$ onto $T_1^*S^3\times\bR_{>0}$.

\begin{thm}
\label{thm:S-maps-to-geodesics}
The map $S$ sends a negative-energy Kepler orbit to a curve in $T^*_1 S^3\times \bR_+$ that projects to a geodesic in $T_1^*S^3$ and a constant in $\bR_+$.
More precisely, if $\phi(t)=(q(t),p(t))$ is an orbit of $X_H$ with $H<0$ and
\[
S(\phi(t))=(r(t),s(t),\sqrt{-2H}),
\]
then there exists a function $\theta(t)$ such that
\[
\begin{cases}
r(t)=\cos\theta(t)\,r(0)+\sin\theta(t)\,s(0),\\
s(t)=-\sin\theta(t)\,r(0)+\cos\theta(t)\,s(0),
\end{cases}
\qquad \dot\theta(t)=\dfrac{\sqrt{-2H}}{|q(t)|}.
\]
\end{thm}

\Proof
Using the Kepler equations $\dot q=p$, $\dot p=-{q}/{|q|^3}$,
one computes the time derivatives of $r_0,\bar r,s_0,\bar s$:

\smallskip
\noindent\emph{(1) Derivative of $r_0$:}
\[
\dot r_0
=|q|\,\frac{d}{dt}|p|^2+|p|^2\,\frac{d}{dt}|q|
=\frac{\bra q,p\ket}{|q|}\,(2H)
=\frac{\sqrt{-2H}}{|q|}\,s_0.
\]

\noindent\emph{(2) Derivative of $\bar r$:}
\[
\dot{\bar r}
=\sqrt{-2H}\Big(\frac{d|q|}{dt}\,p+|q|\,\dot p\Big)
=\frac{\sqrt{-2H}}{|q|}\Big(\bra q,p\ket\, p-\frac{q}{|q|}\Big)
=\frac{\sqrt{-2H}}{|q|}\,\bar s.
\]

\noindent\emph{(3) Derivative of $s_0$:}
\[
\dot s_0
=-\sqrt{-2H}\,\frac{d}{dt}\bra q,p\ket
=-\sqrt{-2H}\Big(|p|^2-\frac1{|q|}\Big)
=-\frac{\sqrt{-2H}}{|q|}\,(|p|^2|q|-1)
=-\frac{\sqrt{-2H}}{|q|}\,r_0.
\]

\noindent\emph{(4) Derivative of $\bar s$:}
\begin{align*}
\dot{\bar s}
&=-\frac{d}{dt}\Big(\frac{q}{|q|}\Big)
+\frac{d}{dt}\big(\bra q,p\ket\, p\big)\\
&=\Big(|p|^2-\frac{2}{|q|}\Big)p
=-\frac{\sqrt{-2H}}{|q|}\,\bar r.
\end{align*}

These four relations are exactly the equations of motion for a geodesic on $S^3$
parametrized by the angle $\theta(t)$ solving 
$\dot\theta(t)=\frac{\sqrt{-2H}}{|q(t)|}$.
\qed

\subsection{Eccentric anomaly and the Ligon--Schaaf rotation}
\label{sec:LS-negative}

Let $\phi(t)$ be a bound Kepler orbit with energy $H<0$ and periapsis at $t=0$.
Denote by $a=-\frac{1}{2H}$ the semi--major axis and by $e$ the eccentricity.
The \emph{eccentric anomaly} $\varepsilon(t)$ and \emph{mean anomaly} $M(t)$ satisfy
\[
|q(t)|=a(1-e\cos\varepsilon(t)),
\qquad
M(t)=\varepsilon(t)-e\sin\varepsilon(t)=\omega t,\quad
\omega=\frac{2\pi}{T}=a^{-3/2}.
\]

\begin{thm}
\label{pro:theta-is-eccentric}
In the situation of Theorem~\ref{thm:S-maps-to-geodesics}, the angle $\theta(t)$
is equal to the eccentric anomaly $\varepsilon(t)$ of the orbit $\phi(t)$.
\end{thm}

\Proof
We have
\[
\dot\theta(t)=\frac{\sqrt{-2H}}{|q(t)|}
=\frac{1}{a^{1/2}|q(t)|}
=\frac{1}{a^{3/2}(1-e\cos\varepsilon(t))}.
\]
On the other hand
\[
M(t)=\varepsilon(t)-e\sin\varepsilon(t)=a^{-3/2}t
\quad\Longrightarrow\quad
\dot M=\dot\varepsilon(1-e\cos\varepsilon)=a^{-3/2},
\]
hence
\[
\dot\varepsilon(t)=\frac{1}{a^{3/2}(1-e\cos\varepsilon(t))}
=\dot\theta(t).
\]
Since we choose $\theta(0)=\varepsilon(0)$ at periapsis, we obtain
$\theta(t)\equiv\varepsilon(t)$.
\qed

Thus, the geodesic angle induced by Moser regularization along a Kepler orbit is exactly the
eccentric anomaly. This is the key observation behind the anomaly interpretation of the
Ligon--Schaaf rotation.

To parametrize all geodesics with constant speed (i.e. mean motion), we introduce
an additional rotation that replaces the eccentric anomaly $\varepsilon$ by the
mean anomaly $M$; this is precisely the Ligon--Schaaf modification.

Define
\[R\colon\T^*_1S^3\times\bR_{>0}\to N_-\]
\[
R(r,s)=
\begin{pmatrix}
x_0&\vec x\\
y_0&\vec y
\end{pmatrix}=
\begin{pmatrix}
\cos(-s_0)&-\sin(-s_0) \\
\frac{\sin(-s_0)}{\sqrt{-2H}}&\frac{\cos(-s_0)}{\sqrt{-2H}}
\end{pmatrix}
\begin{pmatrix}
r_0&\vec r\\
s_0&\vec s
\end{pmatrix},
\]
where $(r,s)$ are as in $S(q,p)$.

\begin{thm}
\label{thm:angle-eps-minus-M}
The rotation angle of $R$ is $\varepsilon-M$, where $\varepsilon$ is the eccentric
anomaly and $M$ is the mean anomaly.
\end{thm}

\Proof
To identify the angle, consider a planar Kepler orbit $q=(x,y)$ with $r=|q|$.
Then
\[
\bra q,p\ket=x\dot x+y\dot y=r\dot r.
\]
Using the parametrization by eccentric anomaly, one has
\[
r=a(1-e\cos\varepsilon),\qquad
\dot r=a e\sin\varepsilon\,\dot\varepsilon
\]
and from the computation in Proposition~\ref{pro:theta-is-eccentric},
\[
\dot\varepsilon=\frac{1}{a^{3/2}(1-e\cos\varepsilon)}.
\]
Hence
\[
\bra q,p\ket=r\dot r
=a(1-e\cos\varepsilon)\cdot
\frac{ae\sin\varepsilon}{a^{3/2}(1-e\cos\varepsilon)}
=a^{1/2}e\sin\varepsilon.
\]
Since $H=-\tfrac{1}{2a}$, we have $\sqrt{-2H}=a^{-1/2}$ and therefore
\[
\sqrt{-2H}\,\bra q,p\ket=e\sin\varepsilon=\varepsilon-M.
\]
This is precisely the rotation angle encoded by $s_0$ in the definition of $R$,
so the rotation introduced by $R$ is by the angle $\varepsilon-M$.
\qed

\subsection{The Ligon--Schaaf regularization on negative energies}

We can now define the global Ligon--Schaaf regularization introduced in \cite{LigonSchaaf1976Kepler} on the negative-energy region
by composing $S$ with $R$.

\begin{dfn}
The \emph{Ligon--Schaaf regularization} on $\Sigma_-$ is
\[
\varphi:=R\circ S\colon \Sigma_-\to N_-.
\]
\end{dfn}

Writing $\Theta:=\sqrt{-2H}\,\bra q,p\ket$, the explicit formula for $\varphi$ is
\begin{align*}
x_0&=\sqrt{-2H}\,\bra q,p\ket\sin\Theta+(|q||p|^2-1)\cos\Theta,\\
\vec x&=\Big(\frac{q}{|q|}-\bra q,p\ket\,p\Big)\sin\Theta+\sqrt{-2H}\,|q|p\cos\Theta,\\
y_0&=-\bra q,p\ket\cos\Theta+\frac{1}{\sqrt{-2H}}\,(|q||p|^2-1)\sin\Theta,\\
\vec y&=\frac{1}{\sqrt{-2H}}\Big(\frac{q}{|q|}-\bra q,p\ket\,p\Big)\cos\Theta+|q|p\sin\Theta.
\end{align*}

\begin{thm}\cite{CushmanBates2015Global}
The map $\varphi$ is a symplectomorphism from $(\Sigma_-,dq_i\wedge dp^i)$ to
$(N_-,dx_i\wedge dy^i)$.
\end{thm}

\Proof
Symplecticity follows from the composition of the symplectomorphisms $S$ and $R$.\qed

\begin{thm}
Along each Kepler orbit, $\varphi$ transforms the Kepler flow into a geodesic 
on $T_1^*S^3$ parametrized by the mean anomaly $M$.
\end{thm}

\Proof
If $\phi(t)=(q(t),p(t))$ is a Kepler orbit with negative energy, then by
Theorem~\ref{thm:S-maps-to-geodesics} and Proposition~\ref{pro:theta-is-eccentric}
we have
\[
S(\phi(t))=(r(t),s(t))
=(\cos\varepsilon(t)\,r(0)+\sin\varepsilon(t)\,s(0),\ -\sin\varepsilon(t)\,r(0)+\cos\varepsilon(t)\,s(0)).
\]
Applying $R$ rotates by $\varepsilon-M$, so that
\[
\varphi(\phi(t))=(x(t),y(t))
=(\cos M(t)\,r(0)+\sin M(t)\,s(0),\ -\sin M(t)\,r(0)+\cos M(t)\,s(0)),
\]
i.e.\ a uniform geodesic motion with angle equal to the mean anomaly $M$.
\qed

\begin{thm}\cite{LigonSchaaf1976Kepler}
Under the Ligon--Schaaf map $\varphi$, the integrals of motion transform as
\[
\frac{A_i}{\sqrt{-2H}}\circ\varphi^{-1}=x_i y_0-x_0 y_i,\qquad
L_i\circ\varphi^{-1}=\varepsilon_{ijk}x_j y_k,\qquad
H\circ\varphi^{-1}=-\frac{1}{2|y|^2}.
\]
\end{thm}

\begin{cor}
\[
A_i^*=\frac{A_i}{\sqrt{-2H}}(q,p)\circ\varphi^{-1}=[A_i-(H+\tfrac{1}{2})q]\circ \MM=A_i',\qquad\]
\[L_i^*=L_i\circ\varphi^{-1}=L_i\circ \MM=L_i',\qquad\]
\[2(-2H(q,p)\circ\varphi^{-1})^{-\frac{1}{2}}-2=K_+(q,p)\circ \MM.
\]
\end{cor}

\begin{tikzcd}
{(N_-,(x,y))} \arrow[r, "R^{-1}"', shift right] \arrow[rr, "\varphi^{-1}", bend left, shift left] 
& \T^*_1S^3\times \bR \arrow[l, "R"', shift right] \arrow[r, "S^{-1}"', shift right] 
& {(\Sigma_-,(q,p))} \arrow[l, "S"', shift right] \arrow[ll, "\varphi", bend left, shift left] \\
\\
\\
{(\T^*(\bR^3\!\setminus\!\{0\}),(q,p))} \arrow[rr, "\MM"', shift right] & & 
{(N_-,(x,y))} \arrow[ll, "\MM^{-1}"', shift right]
& \text{where } \Sigma_-\subsetneq T^*(\bR^3\setminus\{0\}).
\end{tikzcd}

\clearpage

\section{Positive energy: Belbruno and hyperbolic Ligon--Schaaf}
\label{sec:positive}

We now turn to the positive-energy region, following Belbruno's hyperbolic regularization
and extending it in direct analogy with the negative-energy case.

\subsection{Belbruno’s hyperbolic regularization}

We replace the 3-sphere by the hyperbolic space
\[
H^3_+ = \{ (x_0, x_1, x_2, x_3) \in \mathbb{R}^{1,3} 
  \mid x_0^2 - x_1^2 - x_2^2 - x_3^2 = 1,\; x_0 > 0 \}.
\]
which is a subspace of Minkowski space with metric $$|x|^2=\bra x,x\ket_{(1,3)} := x_0^2-x_1^2-x_2^2-x_3^2,$$ the induced metric 
and the metric on the cotangent space is also Minkowski metric, notice that $|y|^2<0$. 
Let $b\colon H_+^3\to\bR^3$ be the stereographic projection adapted to this setting.
As before, a cotangent lift $b_\#$ and the symplectic map $g(q,p)=(-p,q)$ combine into
the \emph{Belbruno map} introduced in \cite{Belbruno1980Regularization} and \cite{Knorrer1997Kepler}
\[
\mathcal{B}:=g\circ b_\# \colon T^*H^3_+ \longrightarrow T^*\mathbb{R}^3.
\]
\[
\mathcal{B}(x,y) = \bigl( -(1-x_0)(y_1,y_2,y_3) - y_0(x_1,x_2,x_3),
\; (1-x_0)^{-1}(x_1,x_2,x_3) \bigr).
\]
\[
\mathcal{B}^{-1}(q,p)=\Big((\frac{|p|^2+1}{|p|^2-1},\frac{2p}{1-|p|^2}),(-\bra q,p\ket, \frac{1-|p|^2}{2}q+\bra q,p\ket p)\Big)
\]
Let $\mathcal{B}(x,y)=(q,p)$. A short computation yields
\[
|q|^2 = -(1-x_0)^2|y|^2,\qquad
|q| = \sqrt{-|y|^2}(x_0-1),
\]
and
\[
|p|^2 - 1 = 2(x_0 - 1)^{-1}.
\]
Hence
\[
K_-(q,p) = |q|(|p|^2-1)-2
=2\sqrt{-|y|^2} - 2,
\]
so that $K_-(q,p)=0$ if and only if $|y|^2=-1$. Thus
\[
\mathcal{B}|_{T^*_1H^3_+}\colon T^*_1H_+^3:=
\bigl\{ (x,y)\in \mathbb{H}^2 \,\big|\,
\bra x,y\ket = 0,\;
\bra x,x\ket = 1,\;
\bra y,y\ket = -1,\;
x_0>0 \bigr\}
\longrightarrow
K_-^{-1}(0)
\]
is a diffeomorphism.

\subsection{Group action and momentum map}
Let  $$G_-:=\{(\alpha,\beta)\in\bH\times\bH\mid|\alpha|^2-|\beta|^2=1,\ \alpha\beta^*+\beta\alpha^*=0\}$$
We identify $$\bR^4\cong\bR\times\bH_{pure}\subseteq\bH\times\bH$$
and $$(\bH\times\bH,*)\cong\bH\otimes\bC\cong M_2(\bC),\ \ (\alpha,\beta)*(\alpha',\beta')=(\alpha\alpha'-\beta\beta',\alpha\beta'+\beta\alpha')$$
by $$(\alpha,\beta)\longleftrightarrow(\alpha_0+\sqrt{-1}\beta_0)I_2-\sqrt{-1}(\alpha_i+\sqrt{-1}\beta_i)\sigma^i $$
where $$\sigma^1=\begin{pmatrix}
    0&1\\
    1&0
\end{pmatrix},\ \sigma^2=\begin{pmatrix}
    0&-i\\
    i&0
\end{pmatrix},\ \sigma^3=\begin{pmatrix}
    1&0\\
    0&-1
\end{pmatrix}.$$ Then $SL_2(\bC)\cong G_-$ and $\bR^4\cong\bR\times\bH_{pure}$ is isomorphic to the set of Hermitian matrices.

Hence, the natural action of $SL_2(\bC)$ on $2\times2$ Hermitian matrices induces an
action on $T^*H_+^3$ and hence $G_-$ act on $T^*\bR^3$ via $\mathcal{B}$ by
\[
(\alpha,\beta)\circ(q,p)
=
\bigl( (\beta p + \alpha)\,q\,(\beta p + \alpha)^*,\;
       (\alpha p - \beta)\,(\beta p + \alpha)^{-1} \bigr)
.\]

\begin{thm}\cite{Knorrer1997Kepler}
The formula above defines a symplectic $G_-$–action on 
$T^*\mathbb{R}^3 \cong \mathbb{H}_{\mathrm{pure}}^2$, which preserves $K_-$.
\end{thm}

\Proof
The proof is parallel to that of Theorem~\ref{thm:Gplus-action}, using the Minkowski
metric in place of the Euclidean one and the Belbruno map $\mathcal{B}$ in place
of the Moser map $\MM$.
\qed

\subsection{Momentum map and integrals of motion}

By Noether’s theorem, the $G_+$–action yields a conserved momentum map.
We compute it explicitly and identify it with the Laplace--Runge--Lenz vector and
the angular momentum.

For $(\alpha, \beta)\in G_-$, write $\alpha=1+a\varepsilon,\ \beta=1+b\varepsilon$, where $(a,b)\in\mathfrak{g}_-$. Then 
$$(1+a\varepsilon)(1+a^*\varepsilon)-b\varepsilon b^*\varepsilon=1, (1+a\varepsilon)b^*\varepsilon+b\varepsilon(1+a^*\varepsilon)=0
\Rightarrow a+a^*=0,\ b+b^*=0$$

\begin{pro}
\label{pro:moment-map-negative}
By the above identification, the Lie algebra $\mathfrak{g}_-=\mathrm{Lie}(G_-)$ is
$$
\mathfrak{g}_-\cong\bH_{\mathrm{pure}}\oplus\bH_{\mathrm{pure}}
$$
and the vector field at $(x,y)\in T^*H_+^3$ associated to $(a,b)\in\mathfrak{g}_-$
is
\[
V^{(a,b)}(x,y)=(-2\Re(bx),2x_0b+2\Im(ax)).
\]
The corresponding momentum map $\mu\colon T^*H_+^3\to\mathfrak{g}_-^*$ satisfies
\[
\mu_{(a,b)}(x,y)
=\bra y,-2\Re(bx)\ket_{(1,3)}=\mu_{(a,b)}=2b\cdot(x_iy_0-x_0y_i)+2a\cdot(\varepsilon_{ijk} x_j y_k) 
\]
where $A$ is the Laplace--Runge--Lenz vector and $L$ is the angular momentum.
\end{pro}

\begin{thm}\cite{Belbruno1980Regularization}
\[
[A_i - (H - \tfrac{1}{2})q_i]\circ \mathcal{B} = x_i y_0 - x_0 y_i, 
\qquad
L_i(q,p)\circ \mathcal{B} = \varepsilon_{ijk} x_j y_k, 
\qquad
K_-(q,p)\circ \mathcal{B} = 2\sqrt{-|y|^2} - 2.
\]
\end{thm}
\Proof
Identical to the computation in the negative energy case, replacing Euclidean with hyperbolic sign conventions.
\qed

\begin{pro}
\label{pro:Kminus-XH}
For
\[
H(q,p) = \frac{1}{2}|p|^2 - \frac{1}{|q|}, 
\qquad 
K_-(q,p) = |q|(|p|^2 - 1) - 2,
\]
one has
\[
\{K_-,H\} = \frac{\bra q,p\ket}{|q|^2}\,K_-(q,p),
\]
so in particular $\{K_-,H\}=0$ on $K_-^{-1}(0)=H^{-1}(\tfrac12)$.
Moreover,
\[
X_{K_-} = 2|q|\, X_H \quad\text{on }K_-^{-1}(0).
\]
\end{pro}

\Proof
The computation is identical to that in Proposition~\ref{pro:Kplus-XH}.
\qed

\subsection{Energy-uniform Belbruno map}

We now define the positive-energy analogue of the map $S$.

Let
\[
\Sigma_+:=\{(q,p)\in T^*(\bR^3\setminus\{0\})\mid H(q,p)>0\}
\]
and let $N_+$ be the corresponding subset of $T^*H_+^3$.
Define
\[
S_+(q,p)
=
\Bigl[(r_0,\bar r),\ (s_0,\bar s),\ \sqrt{2H}\Bigr]
\]
with
\[
r_0=|p|^2|q|-1,\qquad \bar r=-\sqrt{2H}|q|p,
\]
\[
s_0=-\sqrt{2H}\langle q,p\rangle,\qquad
\bar s=-(\frac{q}{|q|}-\langle q,p\rangle p).
\]

\begin{pro}
When $H=\tfrac{1}{2}$ one has $S_+|_{H=\frac{1}{2}}=b_\#\times\{1\}$, so $S_+$ extends Belbruno’s map in an energy-uniform manner. Moreover, $S_+$ is the unique map with this property and is invariant under the scaling in Chapter 2.
\end{pro}

The proof is the positive-energy analogue of the argument for $S$, and we omit it.

\begin{thm}
Let $\nu=\frac{1}{\sqrt{2H}}$ and $\theta=-\nu(\bra s,dr\ket+ds_0$), then $S_+$ is a symplectomorphism from $(\Sigma_+,d(\bra q,dp\ket))$ to $(T^*H^3_+\times\bR_+,d\theta)$.
\end{thm}
\Proof
Note that 
\[
d(|p|^2|q|)=2|q|\bra p,dp\ket+\frac{|p|^2}{|q|}\bra q,dq\ket,\ 
d(\nu^{-1}|q|p)=-\nu^{-2}|q|p\, d\nu+\nu^{-1}p\frac{\bra q,dq\ket}{|q|}+\nu^{-1}|q|dp
\]

Then, 
\begin{align*}
S_+^*(\bra s,dr\ket)&=(-\nu^{-1}\bra q,p\ket,-(\frac{q}{|q|}-\bra q,p\ket \, p)), (d(|p|^2|q|), d(\nu^{-1}|q|p))\\
&=-\nu^{-1}\bra q,dp\ket-\nu^{-2}\bra q,p\ket\, d\nu
\end{align*}
\[
\Rightarrow S_+^*(-\nu\bra s,dr\ket)=\bra q,dp\ket+\nu^{-1}\bra q,p\ket\, d\nu
\]

And
\[
S_+^*(s_0)=\nu^{-1}\bra q,p\ket\Rightarrow S_+^*(\nu ds_0)=\nu d\nu^{-1}\bra q,p\ket=d\bra q,p\ket+\bra q,p\ket(-\nu^{-1})d\nu
\]
\[
\Rightarrow S_+^*\theta=\bra q,dp\ket-d\bra q,p\ket
\]
\[
\Rightarrow S_+^*d\theta=dq_i\wedge dp^i
\]\qed

\begin{thm}
    In the hyperbolic space $\mathbb H^n \subset \mathbb R^{1,n}$,
every geodesic is the intersection of $\mathbb H^n$ with a $2$--dimensional
linear subspace of $\mathbb R^{1,n}$ passing through the origin.
\end{thm}

\Proof
We work in Minkowski space $(\mathbb R^{1,n},\langle\cdot,\cdot\rangle_L)$ with
\[
\langle x,y\rangle_L = x_0y_0 - x_1y_1 - \cdots - x_ny_n,
\]
and define the hyperboloid model by
\[
\mathbb H^n = \{ x\in\mathbb R^{1,n} \mid \langle x,x\rangle_L = 1,\ x_0>0 \}.
\]

\medskip
\noindent
\textbf{Step 1: The plane determined by initial data.}
Let $p\in\mathbb H^n$ and let $v\in T_p\mathbb H^n$ be a unit tangent vector.
Then
\[
\langle p,p\rangle_L = 1, \qquad
\langle v,v\rangle_L = -1, \qquad
\langle p,v\rangle_L = 0.
\]
Define
\[
W := \mathrm{span}\{p,v\} \subset \mathbb R^{1,n}.
\]
The restriction of $\langle\cdot,\cdot\rangle_L$ to $W$ has signature $(-,+)$,
so $W$ is a nondegenerate $2$--dimensional linear subspace.

\medskip
\noindent
\textbf{Step 2: Reflection across the plane $W$.}
Let $W^\perp$ denote the Minkowski orthogonal complement of $W$.
Since $W$ is nondegenerate, we have the orthogonal decomposition
\[
\mathbb R^{1,n} = W \oplus W^\perp .
\]
Every $x\in\mathbb R^{1,n}$ can be written uniquely as $x=w+z$ with
$w\in W$ and $z\in W^\perp$.

Define a linear map $S:\mathbb R^{1,n}\to\mathbb R^{1,n}$ by
\[
S(w+z) := w - z , \qquad w\in W,\ z\in W^\perp .
\]
Geometrically, $S$ is the reflection fixing $W$ pointwise and reversing
the orthogonal directions.

\medskip
\noindent
\textbf{Step 3: $S$ is a Lorentz isometry.}
For $x=w+z$ and $y=w'+z'$ with respect to the above decomposition, orthogonality
gives
\[
\langle x,y\rangle_L = \langle w,w'\rangle_L + \langle z,z'\rangle_L.
\]
Since
\[
Sx = w-z, \qquad Sy = w'-z',
\]
we compute
\[
\langle Sx,Sy\rangle_L
= \langle w-z,\, w'-z'\rangle_L
= \langle w,w'\rangle_L + \langle z,z'\rangle_L
= \langle x,y\rangle_L.
\]
Hence $S\in O(1,n)$ and preserves the Minkowski inner product.

\medskip
\noindent
\textbf{Step 4: $S$ preserves $\mathbb H^n$.}
If $x\in\mathbb H^n$, then $\langle x,x\rangle_L=1$, and by the previous step
\[
\langle Sx,Sx\rangle_L=\langle x,x\rangle_L=1.
\]
Thus $S$ maps the two--sheeted hyperboloid to itself.
Moreover, since $p\in W$, we have $S(p)=p$, which lies in the upper sheet
$x_0>0$. By continuity and connectedness of $\mathbb H^n$, it follows that
\[
S(\mathbb H^n)=\mathbb H^n.
\]
Therefore $S$ restricts to an isometry of $\mathbb H^n$.

\medskip
\noindent
\textbf{Step 5: Fixed point set.}
By construction,
\[
Sx = x \iff z=0 \iff x\in W.
\]
Hence the fixed point set of $S$ in $\mathbb H^n$ is
\[
\mathrm{Fix}(S|_{\mathbb H^n}) = W\cap\mathbb H^n .
\]

\medskip
\noindent
\textbf{Step 6: Uniqueness of geodesics.}
Let $\gamma$ be the geodesic in $\mathbb H^n$ with initial data
\[
\gamma(0)=p, \qquad \gamma'(0)=v.
\]
Since $S$ is an isometry, $S\circ\gamma$ is also a geodesic. Moreover,
\[
(S\circ\gamma)(0)=S(p)=p, \qquad
(S\circ\gamma)'(0)=dS_p(v)=v,
\]
because $S$ acts as the identity on $W$ and $v\in W$.

By uniqueness of solutions to the geodesic initial value problem, we conclude
\[
S\circ\gamma = \gamma.
\]
Thus $\gamma(t)$ is fixed by $S$ for all $t$, and hence
\[
\gamma(\mathbb R) \subset W\cap\mathbb H^n.
\]

Since $W\cap\mathbb H^n$ is a connected $1$--dimensional submanifold, it follows
that $\gamma(\mathbb R)=W\cap\mathbb H^n$.\qed

\begin{thm}
The map $S_+$ sends a positive-energy Kepler orbit to a geodesic in $T^*H_+^3$.
\end{thm}

\Proof
The proof is parallel to that of Theorem~\ref{thm:S-maps-to-geodesics}.
Using $\dot q=p$, $\dot p=-q/|q|^3$ and the definitions of $r_0,\bar r,s_0,\bar s$,
one obtains the system
\[
\dot r_0=-\frac{\sqrt{2H}}{|q|}\,s_0,\qquad
\dot{\bar r}=-\frac{\sqrt{2H}}{|q|}\,\bar s,
\]
\[
\dot s_0=-\frac{\sqrt{2H}}{|q|}\,r_0,\qquad
\dot{\bar s}=-\frac{\sqrt{2H}}{|q|}\,\bar r,
\]
which is the equation for a geodesic on $H^3$ parametrized by an angle
$\theta(t)$ with $\dot\theta(t)=\sqrt{2H}/|q(t)|$.
\qed

\subsection{Hyperbolic anomaly and positive Ligon--Schaaf map}

For a positive-energy orbit, one recalls the \emph{hyperbolic anomaly} $\psi_h$
and mean anomaly $M_h$ via
\[
|q|=|a|(e\cosh\psi_h-1),
\qquad
M_h=e\sinh\psi_h-\psi_h=\omega_h t,\quad
\omega_h=|a|^{-3/2}.
\]

\begin{pro}
Let $S(\phi(t))=(r(t),s(t))$, then 
\[
\begin{cases}
r(t)=\cosh\theta(t)r(0)+\sinh\theta(t)s(0),\\
s(t)=\sinh\theta(t)r(0)+\cosh\theta(t)s(0),
\end{cases}
\qquad \dot\theta(t)=-\frac{\sqrt{2H}}{|q(t)|}.
\]
In particular, $\theta(t)$ is equal to the hyperbolic anomaly $\psi_h(t)$ of $\phi(t)$.
\end{pro}

\Proof
As before,
\[
\dot\theta(t)=-\frac{\sqrt{2H}}{|q(t)|}
=-\frac{1}{|a|^{1/2}|q(t)|}
=-\frac{1}{|a|^{3/2}(e\cosh\psi_h(t)-1)}.
\]
Kepler's law gives $M=|a|^{-3/2}t$, hence
\[
\dot M_h=\dot\psi_h(e\cosh\psi_h-1)=|a|^{-3/2}
\quad\Rightarrow\quad
\dot\psi_h=\frac{1}{|a|^{3/2}(e\cosh\psi_h-1)}=-\dot\theta.
\]
With the choice $\theta(0)=-\psi_h(0)$, we get $\theta(t)\equiv-\psi_h(t)$.
\qed

We now define the hyperbolic analogue of the Ligon--Schaaf rotation.
Set
\[
R_+(r,s)=
\begin{pmatrix}
x_0&\bar x\\
y_0&\bar y
\end{pmatrix}=
\begin{pmatrix}
\cosh(s_0)&-\sinh(s_0) \\
-\frac{\sinh(s_0)}{\sqrt{2H}}&\frac{\cosh(s_0)}{\sqrt{2H}}
\end{pmatrix}
\begin{pmatrix}
r_0&\bar r\\
s_0&\bar s
\end{pmatrix}.
\]

\begin{thm}
$R_+$ is a symplectomorphism from $(T^*H^3_+\times\bR_+,d\theta)$ to $(T^*H_+^3,d(\bra x,dy\ket)$.
\end{thm}
\Proof
Note that $R_+^*(y)=\nu(-\sinh(s_0)r+\cosh(s_0)s)$
\[
\Rightarrow d(R^*y)=\nu(-\cosh(s_0)(ds_0)r-\sinh(s_0)dr+\sinh(s_0)(ds_0)s+\cosh(s_0)ds)+(-\sinh(s_0)r+\cosh(s_0)s)d\nu
\]
Therefore, 
\[
R_+^*(\bra x,dy\ket)=\bra R_+^*x,d(R_+^*y)\ket=-\nu(\bra s,dr\ket+ds_0)
\]
\[
\Rightarrow R_+^*(dx_i\wedge dy^i)=d\theta
\]\qed

\begin{thm}
The rotation parameter of $R_+$ is
\[
-\psi_h-M_h,
\]
where $\psi_h$ is the hyperbolic anomaly and $M_h$ the mean anomaly.
\end{thm}

\Proof
To compute the angle, let $q=(x,y)$ be a planar Kepler orbit, $r=|q|$.
Then
\[
\bra q,p\ket=x\dot x+y\dot y=r\dot r
=|a|(e\cosh\psi_h-1)\cdot\frac{e\sinh\psi_h}{|a|^{1/2}(e\cosh\psi_h-1)}
=|a|^{1/2}e\sinh\psi_h.
\]
Since $H=\frac{1}{2|a|}$, one has $\sqrt{2H}=|a|^{-1/2}$ and thus
\[
-\sqrt{2H}\,\bra q,p\ket=-e\sinh\psi_h=-M_h-\psi_h.
\]
This determines the rotation parameter occurring in $R_+$.
\qed

\begin{dfn}
The \emph{positive-energy Ligon--Schaaf map} \cite{LigonSchaaf1976Kepler} is the composition
\[
\varphi_+:=R_+\circ S_+\colon\Sigma_+\to N_+.
\]
\end{dfn}

Writing $\Theta:=-\sqrt{2H}\bra q,p\ket$, one obtains the explicit formula
\begin{align*}
x_0&=-\sqrt{2H}\bra q,p\ket\sinh\Theta+(|q||p|^2-1)\cosh\Theta,\\
\bar x&=\Big(\frac{q}{|q|}-\bra q,p\ket \, p\Big)\sinh\Theta-\sqrt{2H}|q|p\cosh\Theta,\\
y_0&=-\bra q,p\ket\cosh\Theta-\frac{1}{\sqrt{2H}}(|q||p|^2-1)\sinh\Theta,\\
\bar y&=\frac{-1}{\sqrt{2H}}\Big(\frac{q}{|q|}-\bra q,p\ket \, p\Big)\cosh\Theta+|q|p\sinh\Theta.
\end{align*}

\begin{thm}
The map $\varphi_+$ is a symplectomorphism from $(\Sigma_+,dq_i\wedge dp^i)$
to $(N_+,dx_i\wedge dy^i)$.
If $\varphi_+(q(t),p(t))=(x(t),y(t))$ for a Kepler orbit $(q(t),p(t))$, then
\[
\begin{cases}
x(t)=\cosh(-\psi_h(t)-s_0(t))r(0)+\sinh(-\psi_h(t)-s_0(t))s(0),\\
y(t)=\dfrac{-\sinh(-\psi_h(t)-s_0(t))r(0)+\cosh(-\psi_h(t)-s_0(t))s(0)}{\sqrt{2H}},
\end{cases}
\]
and $-\psi_h-s_0=-\psi_h+(\psi_h+M_h)=M_h$.
Thus, the geodesic motion is parametrized by the mean anomaly $M_h$.
\end{thm}

\Proof
The first claim follows as in the negative-energy case, since both $S_+$ and $R_+$
are symplectic. The formulas for $(x(t),y(t))$ are obtained by composing the
hyperbolic rotation with parameter $\psi_h$ from $S_+$ with the additional rotation
encoded by $R_+$. The identity $\psi_h-s_0=-M_h$ follows from the previous theorem.
\qed

As before, one shows that under $\varphi_+$ the Laplace--Runge--Lenz vector and
angular momentum transform into linear functions in $(x,y)$, and that
$H\circ\varphi_+^{-1}=\frac{1}{2|y|^2}$.

\begin{thm}\cite{LigonSchaaf1976Kepler} 
\[
\frac{A_i}{\sqrt{2H}}(q,p)\circ\varphi_+^{-1}=x_0y_i-x_iy_0,\qquad\]
\[L_i(q,p)\circ\varphi_+^{-1}=\varepsilon_{ijk}x_iy_j,\qquad\]
\[H(q,p)\circ\varphi_+^{-1}=\frac{1}{2|y|^2}.
\]
\end{thm}

\begin{cor}
\[
A_i^*=\frac{A_i}{\sqrt{2H}}(q,p)\circ\varphi_+^{-1}=[A_i-(H-\tfrac{1}{2})q]\circ \mathcal{B}=A_i',\qquad\]
\[L_i^*=L_i\circ\varphi_+^{-1}=L_i\circ \mathcal{B}=L_i',\qquad\]
\[2\sqrt{(2H(q,p)\circ\varphi_+^{-1})^{-1}}-2=K_-(q,p)\circ \mathcal{B}.
\]
\end{cor}

\begin{tikzcd}
{(N_+,(x,y))} \arrow[r, "R_+^{-1}"', shift right] \arrow[rr, "\varphi_+^{-1}", bend left, shift left] 
& \T^*_1H^3\times \bR \arrow[l, "R_+"', shift right] \arrow[r, "S_+^{-1}"', shift right] 
& {(\Sigma_+,(q,p))} \arrow[l, "S_+ "', shift right] \arrow[ll, "\varphi_+", bend left, shift left] \\
\\
\\
{(\T^*(\bR^3\!\setminus\!\{0\}),(q,p))} \arrow[rr, "\mathcal{B}"', shift right] & & {(N_+,(x,y))} \arrow[ll, "\mathcal{B}^{-1}"', shift right]
& {\rm where\ }\  \Sigma_+\subsetneq\T^*(\bR^3\!\setminus\!\{0\})
\end{tikzcd}

\clearpage
\section{Zero energy and parabolic regularization}
\label{sec:zero}

Finally, we consider the zero--energy level $H=0$, which separates the
negative--energy (elliptic) and positive--energy (hyperbolic) regimes and
corresponds to parabolic Keplerian motion.
From the dynamical point of view, this case may be regarded as a limiting
situation of both the elliptic and hyperbolic problems.

However, the Ligon--Schaaf regularizations constructed in the negative and
positive energy regions become degenerate as $H\to 0^\pm$.
More precisely, if one formally takes the limit of the Ligon--Schaaf coordinates
$(x_0,\bar x,y_0,\bar y)$ as $H\to 0$, one obtains
\[
\lim_{H\to 0^+}\bigl(x_0,\bar x,y_0,\bar y\bigr)
=
\left(
0,\;
-\Bigl(|q|\,p+\langle q,p\rangle\Bigl(\tfrac{q}{|q|}-\langle q,p\rangle\,p\Bigr)\Bigr),\;
0,\;
\tfrac{q}{|q|}-\langle q,p\rangle\,p
\right),
\]
which is independent of time along the orbit.
In particular, the image collapses to a fixed point rather than producing a
nontrivial geodesic motion.

This degeneration indicates that, unlike the elliptic and hyperbolic cases,
the zero--energy level cannot be treated by a direct limit of the
Ligon--Schaaf construction.
Instead, a separate regularization adapted to the parabolic geometry is
required.
In what follows, we introduce a Euclidean regularization which plays the role
of a parabolic analogue of the Moser and Belbruno maps and allows us to recover
a uniform geometric description of parabolic Kepler orbits.

\subsection{A Euclidean regularization}

Similarly to the previous two cases, consider the cotangent lift of $\frac1x$
\[
\mathcal{L}\colon T^*\bR^3\to\bR^3\times\bR^3\cong T^*\bR^3,\qquad
\mathcal{L}(x,y)=(\frac{xyx}{2},2x^{-1}),
\]
with inverse
\[
\mathcal{L}^{-1}(q,p)=(2p^{-1},\frac{pqp}{2})=(-|q|p,\frac{q}{|q|}-\bra q,p\ket p).
\]
In quaternionic notation, $x$ and $y$ are purely imaginary and $xyx$ is again
purely imaginary, so that $\mathcal{L}$ maps $T^*\bR^3$ to $T^*\bR^3$.

We introduce a group
\[
G_0=\{(\alpha,c)\mid\alpha\in\bH_{\text{unit}},\ c^*\alpha+\alpha^*c=0\},
\]
with multiplication
\[
(\alpha',c')(\alpha,c)=(\alpha'\alpha, \alpha'c+c'\alpha).
\]
We let $G_0$ act on $T\bR^3\cong\bR^3\times\bR^3$ by
\[
(\alpha,c)*(x,y)=(c\alpha^{-1}+\alpha x\alpha^{-1},\alpha y\alpha^{-1}).
\]
Transporting this action via $\mathcal{L}$ gives an action on $T^*\bR^3$:
\[
(\alpha,c)\cdot(q,p)=((cp+\alpha)q(cp+\alpha)^*,\alpha p(cp+\alpha)^{-1}).
\]

\begin{thm}
The $G_0$–action on $T^*\bR^3$ is symplectic and preserves
\[
K_0(q,p)=|q||p|^2-2.
\]
Moreover,
\[
(\alpha,c)\cdot \mathcal{L}(x,y)=\mathcal{L}((\alpha,c)*(x,y))
\]
for all $(\alpha,c)\in G_0$ and $(x,y)\in T\bR^3$.
\end{thm}

\Proof
The fact that the actions are compatible via $\mathcal{L}$ follows from a direct
calculation with the quaternionic formulas. Symplecticity and the invariance
of $K_0$ are checked as in the negative and positive energy cases; we omit the details.
\qed

\subsection{Momentum map}

As before, let $\alpha=1+a\varepsilon, c=b\varepsilon$, where $(a,b)\in\mathrm{Lie}(G_0)$, then 
$$(1+a\varepsilon)^2=1,b^*\varepsilon(1+a\varepsilon)+(1+a^*\varepsilon)b\varepsilon=0$$
which implies $a,b\in\bH_{pure}$.
Thus,
\[
\mathrm{Lie}(G_0)=\bH_{\text{pure}}\oplus\bH_{\text{pure}},\ 
\Bigg[\begin{pmatrix}
    a_1\\b_1
\end{pmatrix},
\begin{pmatrix}
    a_2\\b_2
\end{pmatrix}\Bigg]
=\begin{pmatrix}
    a_1a_2-a_2a_1\\
    a_2b_1-a_1b_2+b_2a_1-b_1a_2
\end{pmatrix}
\]
and an infinitesimal action
\[
V^{(a,b)} = (ax-xa+b,\ ay-ya)
\]
on $T\bR^3$. The corresponding moment map is
\[
\mu_{(a,b)} = \iota_V\theta =y\cdot V^i=a\cdot(x\wedge y)+b\cdot y.
\]
In $(q,p)$-coordinates, after composing with $\mathcal{L}$, one finds
\[
\mu_{(a,b)} = b\cdot A+a\cdot L,
\]
where $A$ and $L$ are again the Laplace--Runge--Lenz vector and the angular momentum.

\begin{thm}
On $T\bR^3$ one has
\[
A \circ \mathcal{L}= y, 
\qquad
L\circ \mathcal{L}= x\wedge y, 
\qquad
K_0\circ \mathcal{L} = 2|y| - 2.
\]
\end{thm}

\Proof
$$y=\frac12qpq=\frac12\Im(pq\cdot p)=\frac12(\Re(pq)p+\Im(pa)\wedge p)=\frac12(-\bra p, q\ket p+(p\wedge q)\wedge p)$$
$$=\frac12(-\bra q, p\ket p-p\bra q,p \ket+|p|^2q)=\frac12(-2\bra q,p\ket p+|p|^2q)=\frac{q}{|q|}-\bra q,p\ket p$$
and since $H=0\Leftrightarrow\frac{1}{2}|p|^2|q|=1$,
$$A=p\wedge q\wedge p-\frac{q}{|q|}=-\bra q,p\ket p+\frac{q}{|q|}=y$$
The other two equations are also proved by direct computation.
\qed

\subsection{Poisson bracket}

The analogue of Propositions~\ref{pro:Kplus-XH} and~\ref{pro:Kminus-XH} holds.

\begin{pro}
For
\[
H(q,p) = \frac{1}{2}|p|^2 - \frac{1}{|q|}, 
\qquad 
K_0(q,p) = |q||p|^2 - 2,
\]
one has
\[
\{K_0,H\}=\frac{\bra q,p\ket}{|q|^2}K_0(q,p),
\]
so in particular $\{K_0,H\}=0$ on $K_0^{-1}(0)=H^{-1}(0)$.
Moreover,
\[
X_{K_0}=2|q|\,X_H\quad\text{on }K_0^{-1}(0).
\]
\end{pro}

\Proof
The computation is again identical in spirit to the one for $K_+$ and $K_-$:
\[
\frac{\partial K_0}{\partial q_i} = |p|^2\frac{q_i}{|q|},\quad
\frac{\partial K_0}{\partial p_i} = 2|q|p_i,
\]
so
\[
\{K_0,H\} 
=\frac{|p|^2}{|q|}\bra q,p\ket-2|q|\frac{\bra q,p\ket}{|q|^3}
=\frac{\bra q,p\ket}{|q|^2}K_0(q,p).
\]
On $K_0^{-1}(0)$ one has $|q||p|^2=2$, hence
\[
X_{K_0}
=(2|q|p,-\tfrac{|p|^2}{|q|}q)
=2|q|\Big(p,-\frac{q}{|q|^3}\Big)=2|q|X_H.
\]
\qed

Beacuse the case of zero energy is only one energy shell, we don't need to uniform it. However, we still can get a parallel result to previous two cases.

\subsection{Geodesics and parabolic anomaly}
\label{subsec:zero-geodesic-anomaly}

\begin{thm}
\label{thm:L-maps-to-lines-parallel}
The map $\mathcal{L}$ sends each zero-energy Kepler orbit to a geodesic in $\bR^3$.
More precisely, if $\phi(t)=(q(t),p(t))$ is a Kepler orbit with $H(\phi(t))=0$ and we set
\[
(x(t),y(t)):=\mathcal{L}^{-1}(q(t),p(t))=(2p(t)^{-1},\,\frac{p(t)q(t)p(t)}{2}),
\]
then there exists a scalar function $\theta(t)$ such that
\[
y(t)\equiv y(0),\qquad
x(t)=x(0)+\theta(t)\,y(0),
\qquad
\dot\theta(t)=\frac{1}{2|q(t)|}.
\]
In particular, the curve $x(\theta)$ is an affine line with constant speed,
hence a geodesic in $\bR^3$.
\end{thm}

\Proof
Differentiate $2y=pqp$ with respect to $t$:
\[
2\dot y=\dot p\,q\,p+p\,\dot q\,p+p\,q\,\dot p.
\]
Using $\dot q=p$ and $\dot p=-\frac{q}{|q|^3}$, and the fact that $q$ is purely imaginary
so that $qq=-|q|^2$, we obtain
\[
\dot p\,q\,p=-\frac{q}{|q|^3}q\,p=\frac{1}{|q|}p,\qquad
p\,q\,\dot p=pq\Bigl(-\frac{q}{|q|^3}\Bigr)=\frac{1}{|q|}p,
\]
and since $p$ is purely imaginary, $p^2=-|p|^2$ is real, hence
\[
p\,\dot q\,p=p\,p\,p=p^3=-|p|^2p.
\]
Therefore
\[
2\dot y=\Bigl(\frac{2}{|q|}-|p|^2\Bigr)p.
\]
On $H=0$ we have $|p|^2=\frac{2}{|q|}$, so $\dot y=0$, i.e.\ $y(t)\equiv y(0)$.

Next, $x=2p^{-1}$ satisfies the inverse-derivative identity
\[
\dot x=\frac{d}{dt}(2p^{-1})=-2p^{-1}\dot p\,p^{-1}
=\frac{2}{|q|^3}p^{-1}q\,p^{-1}.
\]
Since $p$ is purely imaginary, $p^{-1}=\bar p/|p|^2=-p/|p|^2$, hence
\[
p^{-1}q\,p^{-1}=\frac{pqp}{|p|^4}=\frac{y}{|p|^4}.
\]
Thus
\[
\dot x=\frac{2}{|q|^3}\frac{y}{|p|^4}.
\]
On $H=0$, $|p|^2=\frac{2}{|q|}$ so $|p|^4=\frac{4}{|q|^2}$, and therefore
\[
\dot x=\frac{2}{|q|^3}\cdot \frac{|q|^2}{4}\,y=\frac{1}{2|q|}\,y.
\]
Since $y$ is constant, there exists a scalar function $\theta(t)$ with
\[
\dot\theta(t)=\frac{1}{2|q(t)|},\qquad
\dot x(t)=\dot\theta(t)\,y(0)
\quad\Longrightarrow\quad
x(t)=x(0)+\theta(t)\,y(0),
\]
which shows that $x(\theta)$ is an affine line with constant direction and speed.
\qed

\medskip
\noindent
We now identify the parameter $\theta(t)$ with the parabolic anomaly.

\begin{thm}
\label{pro:theta-is-parabolic}
Let $\psi_p=\tan(\theta/2)$, where $\theta$ is the true anomaly, be the parabolic anomaly of the planar Kepler motion
(as in Section~3), and let the parabolic mean anomaly be
\[
M_p=\omega_p(t-t_p)=\psi_p+\frac{1}{3}\psi_p^3.
\]
Then the geodesic parameter $\theta(t)$ in Theorem~\ref{thm:L-maps-to-lines-parallel}
satisfies
\[
\theta(t)=\frac{\psi_p}{\omega_p\mathfrak p}
\]
where $\omega_p$ and $\mathfrak{p}$ are constant.
\end{thm}

\Proof
For a parabolic Kepler orbit, one has (with $\mathfrak p$ fixed)
\[
r=\frac{\mathfrak p}{2}(1+\psi_p^2).
\]
By Theorem~\ref{thm:L-maps-to-lines-parallel},
\[
\dot\theta(t)=\frac{1}{2|q(t)|}=\frac{1}{2r(t)}
=\frac{1}{\mathfrak p(1+\psi_p(t)^2)}.
\]
On the other hand, the defining relation
\[
M_p=\psi_p+\frac{1}{3}\psi_p^3
\]
gives
\[
\dot M_p = \dot\psi_p(1+\psi_p^2)=\omega_p.
\Rightarrow
1+\psi_p^2=\frac{\omega_p}{\dot\psi_p}
\]
Therefore,
\[\dot\theta(t)=\frac{\dot\psi_p}{\omega_p\mathfrak p}\]
With the normalization $\theta(t_p)=\psi_p(t_p)=0$, we conclude $\theta(t)=\frac{\psi_p}{\omega_p\mathfrak p}$.
\qed






\section*{Acknowledgements}

I would like to express my sincere gratitude to all those who have provided support and assistance throughout the completion of this thesis.

First and foremost, I would like to thank my advisor, Professor Siye Wu, for his invaluable guidance, patience, and insightful suggestions during the course of this research. His expertise and continuous support were essential to the development and completion of this thesis.

I would also like to extend my sincere appreciation to Professor Kuo-Chang Chen for the knowledge and perspectives gained from his courses, which have been of great benefit to this research. In addition, I am grateful to Professor Nan-Kuo Ho for her valuable comments and guidance throughout the research process.

Finally, I would like to express my deepest gratitude to my family and friends for their constant encouragement and support, which have been a vital source of motivation throughout my graduate studies.

\printbibliography

@article{Moser1970Kepler,
  author  = {Moser, J{\"u}rgen},
  title   = {Regularization of Kepler’s Problem and the Averaging Method on a Manifold},
  journal = {Communications on Pure and Applied Mathematics},
  volume  = {23},
  year    = {1970},
  pages   = {609--636},
}

@article{Belbruno1980Regularization,
  author  = {Belbruno, Edward A.},
  title   = {A New Regularization of the Restricted Three-Body Problem and an Application},
  journal = {Celestial Mechanics},
  volume  = {21},
  year    = {1980},
  pages   = {149--168},
}

@book{CushmanBates2015Global,
  author    = {Cushman, Richard H. and Bates, Larry M.},
  title     = {Global Aspects of Classical Integrable Systems},
  edition   = {2},
  publisher = {Birkh{\"a}user},
  address   = {Basel},
  year      = {2015},
}

@article{LigonSchaaf1976Kepler,
  author  = {Ligon, Thomas and Schaaf, Manfred},
  title   = {On the Global Symmetry of the Classical Kepler Problem},
  journal = {Reports on Mathematical Physics},
  volume  = {9},
  number  = {2},
  year    = {1976},
  pages   = {281--300},
}

@incollection{Knorrer1997Kepler,
  author    = {Kn{\"o}rrer, Horst},
  title     = {Explicit Symmetries of the Kepler Hamiltonian},
  booktitle = {Global Aspects of Classical Integrable Systems},
  editor    = {Cushman, Richard H. and Bates, Larry M.},
  publisher = {Birkh{\"a}user},
  address   = {Basel},
  year      = {1997}
}

@article{KustaanheimoStiefel1965KS,
  author  = {Kustaanheimo, Paul and Stiefel, Eduard},
  title   = {Perturbation Theory of Kepler Motion Based on Spinor Regularization},
  journal = {Journal f{\"u}r die reine und angewandte Mathematik},
  volume  = {218},
  year    = {1965},
  pages   = {204--219}
}

@article{Arnold1997MathematicalAO,
  title={Mathematical aspects of classical and celestial mechanics},
  author={Vladimir I. Arnold and Valery V. Kozlov and Anatoly I. Neishtadt},
  year={1997},
  url={https://api.semanticscholar.org/CorpusID:117778987}
}

\end{document}